\documentclass[12pt]{article}
\usepackage{amssymb,amsmath}
\def\hang{\hangindent\parindent}
\def\tex#1{\indent\llap{[#1]\enspace}\ignorespaces}
\def\re{\par\hang\tex}
\def\Res{{\rm Res}_t}
\def\F{{\mathbb{F}}}
\def\SUM#1#2{\mbox{$\sum\limits_{#1}^{#2}$}}
\def\SUM#1#2{\sum_{#1}^{#2}}
\def\OPLUS#1{\raisebox{-4pt}{\mbox{${\begin{array}{c}\mbox{\Large $\oplus$}\\[-6pt]\sc
#1\end{array}}$}}}
\def\OPLUS#1{\bigoplus_{#1}}
\def\VS#1{}
\def\vs{\vspace*}
\def\a{\alpha}

\def\b{\beta}
\def\d{\delta}
\def\D{\Delta}
\def\BB{\mathcal B}

\def\WW{{\mathcal W}}
\def\LL{{\mathcal B}}

\def\G{\Gamma}

\def\L{\Lambda}

\def\sc{\scriptstyle}
\def\ssc{\scriptscriptstyle}

\def\cl{\centerline}

\def\Rar{\Rightarrow}

\def\AA{{\mathcal A}}

\def\PP{{\mathcal P}}

\def\ni{\noindent}
\def\ptl{\partial}

\def\N{\mathbb{N}{\ssc\,}}
\def\Z{\mathbb{Z}{\ssc\,}}

\numberwithin{equation}{section}
\newtheorem{theo}{Theorem}[section]

\newtheorem{prop}[theo]{Proposition}

\begin{document}
\cl {{\Large\bf Verma modules over a Block Lie algebra}} \vs{6pt}
\cl{Qifen Jiang$^{\,1)}$,\ \ \ Yuezhu Wu$^{\,1,\,2)}$} \cl{\small
$^{1)}$Department of Mathematics, Shanghai Jiao Tong University,
 Shanghai 200240, China}

\cl{\small $^{2)}$Department of Mathematics, Qufu Normal
University,
 Qufu 273165, China}
\vs{10pt}\par

{\small
\parskip .01 truein
\baselineskip 3pt \lineskip 3pt

\noindent{{\bf Abstract.} Let $\BB$ be the Lie algebra with basis
$\{L_{i,j},\,C\,|\,i,j\in\Z\}$ and relations
$[L_{i,j},L_{k,l}]=((j+1)k-i(l+1))L_{i+k,j+l}+i\d_{i,-k}\d_{j+l,-2}C,\,
[C,L_{i,j}]=0.$ It is proved that an irreducible highest weight
$\BB$-module is quasifinite if and only if it is a proper quotient
of a Verma module. For an additive subgroup $\G$ of $\F$, there
corresponds to a Lie algebra $\BB(\G)$ of Block type. Given a total
order $\succ$ on $\G$ and a weight $\L$, a Verma $\BB(\G)$-module
$M(\L,\succ)$ is defined. The irreducibility of $M(\L,\succ)$ is
completely determined. \vs{5pt}

\noindent{\bf Key words:} Verma modules, Lie algebras of Block type,
irreducbility, quasifinite modules}

\noindent{\it Mathematics Subject Classification (2000):} 17B10,
17B65,  17B68.}
\parskip .001 truein\baselineskip 6pt \lineskip 6pt

\vs{6pt}

\par

\cl{\bf\S1. \
Introduction}\setcounter{section}{1}\setcounter{equation}{0}

Since a class of infinite dimensional simple Lie algebras was
introduced by Block [B], generalizations of Block algebras have been
studied by some authors, partially because they are closely related
to the Virasoro algebra or algebras associated quantum plane (e.g.,
[DZ, LW, OZ, S1--S4, SZ, X1, X2, ZM]).

For an additive subgroup $\G$ of a field $\F$ of characteristic $0$,
we consider in this paper the {\it Lie algebra $\BB(\G)$ of Block
type} with basis $\{L_{\a,i},C\,|\,\a \in \G,\, i\in\Z \}$ and
\VS{-5pt}relations
\begin{equation}\label{def}
\begin{array}{ll}
[L_{\a,i},L_{\b,j}]  =
((i+1)\b-(j+1)\a)L_{\a+\b,i+j}+\a\d_{\a,-\b}\d_{i+j,-2}C, \\[4pt]
[C,L_{\a,i}] = 0\VS{-5pt}.
\end{array}
\end{equation}
This particular Block type Lie algebra attracts our attention
because it is also a special case of Cartan type $S$ Lie algebras
(e.g., [X, SX1]), or Poisson (or Hamiltonian) algebras (e.g., [X,
S6, SX2]). One of our motivations to study representations of this
Lie algebra is to better understand representations of Lie algebras
of 4 families of Cartan type.

The Lie algebra $\BB(\G)$ is $\G$-graded (but not finitely
graded)\VS{-7pt}:
 \begin{equation}\label{grad}
 \BB(\G)=\OPLUS{\a\in\G}\BB(\G)_\a,\ \ \mbox{ where \ }
 \BB(\G)_\a={\rm span}\{L_{\a,i}|~i\in \Z\}+\d_{\a,0}\F C\VS{-7pt}.
\end{equation}
Note that $\BB(\G)_0$ is an infinite dimensional commutative
subalgebra of $\BB(\G)$ (but it is not a Cartan subalgebra). For a
total order `` $\succ $'' on $\G$ compatible with its group
structure (in case $\G=\Z$, we always choose the normal order on
$\Z$), we denote $\G_{\pm}=\{x\in\G\,|\,\pm x\succ0\}$. Then
$\G=\G_{+}\cup\bigl\{0\bigr\}\cup\G_{-}$, and we have the {\it
triangular \VS{-5pt}decomposition}
$$\BB(\G)=\BB(\G)_-\oplus \BB(\G)_0 \oplus \BB(\G)_+,
\mbox{ \ where \ } \BB(\G)_{\pm}=\OPLUS{\pm\a\succ
0}\BB(\G)_\a\VS{-5pt}.$$ Thus we have the notion of a Verma module
$M(\L,\succ)$ with respect to a function, called a {\it weight},
$\L\in\BB(\G)_0^*$ (the dual space of $\BB(\G)_0$) and the order
``$\succ$'' (cf.~(\ref{verma-module})).

A weight $\L$ is described by the {\it central charge} $c=\L(C)$ and
its {\it labels} $\L_i=\L(L_{0,i-1})$ for $i\in\Z$. For $j\in\Z$, we
introduce the {\it $j$-th generating \VS{-5pt}series}
\begin{equation}\label{generating-series}
\D_\L^{(j)}(z)=\SUM{i=0}\infty\frac{\L_{i+j}}{i!}z^{i}\VS{-5pt}.
\end{equation}
A function $\D(z)$ is called a {\it quasipolynomial} if it is a
linear combination of functions of the form $p(z)e^{\a z}$, where
$p(z)\in\F[z],\,\a\in\F$. Recall [KL, KR] the well-known
characterization that a formal power series is a quasipolynomial if
and only if it satisfies a nontrivial linear differential equation
with constant coefficients.

A $\BB(\G)$-module $V$ is {\it quasifinite} if $V$ is finitely
$\G$-graded, namely\VS{-5pt}, $$V =\OPLUS{\a\in \G}V_\a\mbox{ \ \
with \ \ } \BB(\G)_\a V_\b\subset V_{\a+\b},\ \ {\rm dim}V_\a<\infty
\mbox{\ \ \ for \ \ }\a,\b\in \G\VS{-7pt}.$$

Quasifinite modules are closely studied by some authors, e.g., [KL,
KR, KWY, LZ, S3, S4]. Quasifinite irreducible modules over the Lie
algebras $\BB'(\G)$ are proved to be a highest or lowest weight
module in [S5], where \begin{equation}\label{b'}\BB'(\G)={\rm
span}\{L_{\a,i},C\,|\,\a \in \G,\, i\geq -1 \}\end{equation} is a
subalgebra of $\BB(\G)$.

Since each grading space $\BB(\G)_\a$ in (\ref{grad}) is
infinite-dimensional, the classification of quasifinite modules is a
nontrivial problem as pointed in [KL, KR]. Thus the classification
of graded modules with infinite dimensional grading spaces is a
nontrivial problem as well.

 The main results
in this paper are the following.

\begin{theo}
\label{main-theo1} Suppose $V$ is an irreducible highest weight
$\BB(\Z)$-module of weight $\L$. The following statements are
equivalent:
\begin{itemize}\parskip-4pt
\item[{\rm(1)}]
$V$ is quasifinite;
\item[{\rm(2)}]
$V$ is a proper quotient of the Verma module $M(\L,\succ)$;
\item[{\rm(3)}]
$\Sigma_\L^{(j)}(z):=z\D_\L^{(-j)}(z)-j\D_\L^{(-j-1)}(z)-\frac{c}{j!}z^j$
is a quasipolynomial and satisfies the same nontrivial linear
differential equation with constant coefficients for all $j\in\Z$.
\end{itemize}
\end{theo}

\begin{theo}
\label{main-theo2} Let $M(\L,\succ)$ be the Verma $\BB(\G)$-module
with highest weight $\L$ $($cf.~$(\ref{verma-module}))$.
\begin{itemize}\parskip-4pt
\item[{\rm(1)}]
With respect to a dense order $``\succ$'' of $\G$
$($cf.~$(\ref{dense}))$, $$\mbox{$M(\L,\succ)$ is irreducible \ \ \
$\Longleftrightarrow$ \ \ \ $\L\neq 0$.}$$Moreover, in case $\L=0$,
the submodule
$$\begin{array}{ll}
M'(0,\succ)=\SUM{k>0,\, \a_1,\cdots,\a_k\in G_+}{}\F
L_{-\a_1,i_1}\cdots L_{-\a_k,i_k}v_0\mbox{ \ \ is irreducible} \\[12pt]
\Longleftrightarrow\ \ \ \forall\,x,y\in \G_+,\
\exists\,n\in\Z_+\mbox{ \ such that \ }nx\succ y.\end{array}$$
 \item[{\rm(2)}] With respect to a
discrete order $`` \succ$" $($cf.~$(\ref{discrete}))$,
$$\mbox{$M(\L,\succ)$ is irreducible}\ \
\Longleftrightarrow\  \ \mbox{$\BB(a\Z)$-module $M_a(\L,\succ)$ is
irreducible}.$$
\end{itemize}
\end{theo}
\vskip10pt

\cl{\bf\S2. \ Proof of the main result
}\setcounter{section}{2}\setcounter{theo}{0}\setcounter{equation}{0}

Let $\G$ be any additive subgroup of $\F$. Denote by $U=U(\BB(\G))$
the universal enveloping algebra of $\BB(\G)$. For any $\L\in
\BB(\G)_0^*$, let $I(\L,\succ)$ be the left ideal of $U$ generated
by the elements
$$\{L_{\a,i}\,|~\a\succ 0, i\in \Z\}\cup\{h-\L(h)\cdot 1\,|~h\in \BB(\G)_0\}.$$
Then the {\it Verma $\BB(\G)$-module} with respect to the order ``
$\succ$''
  is defined as
  \begin{equation}\label{verma-module}M(\L,\succ)=U/I(\L,\succ),\end{equation}
which has  a basis consisting of all vectors of the form
\begin{equation}\label{vector-form}
\begin{array}{lll}
L_{-\a_1,i_1}L_{-\a_2,i_2}\cdots L_{-\a_k,i_k}v_{\L}, &\mbox{where}&
i_j\in \Z,\ 0\prec\a_1\preccurlyeq \cdots\preccurlyeq \a_k,\\[4pt]&\mbox{and}
&i_s\leq i_{s+1}\mbox{ if }\a_s=\a_{s+1}.
\end{array}
\end{equation}
where $v_\L$ is the coset of $1$ in $M(\L,\succ)$. Thus
$M(\L,\succ)$ is a {\it highest weight $\BB(\G)$-module} in the
sense
\begin{equation}\label{highest}M(\L,\succ)=\oplus_{\a\preccurlyeq
0}M_{\a},\end{equation}
 where $M_0=\F v_\L,$ and $M_\a$ for $\a\prec
0$ is spanned by vectors in (\ref{vector-form}) with $\a_1+\cdots+
\a_k=-\a.$ So it is a $\G$-graded $\BB(\G)$-module with
\begin{equation}\label{infty}
{\rm dim}M_{-\a}=\infty\mbox{ \ \ for \ \ }\a \in
\G_+.\end{equation}

 For any $a\in \G$, we  denote
 \begin{equation}\label{BBaZ}\BB(a\Z)={\rm
span}\{L_{na,k},C\,|~n,k\in \Z\},\end{equation} which is a
subalgebra of $\BB(\G)$ isomorphic to $\BB(\Z)$. We also denote
\begin{equation}\label{BBaM}M_a(\L,\succ)\ \ =\ \ \mbox{ $\BB(a\Z)$-submodule of $M(\L,\succ)$
 generated by $v_\L$.}\end{equation}
Denote
\begin{equation}\label{denote-Ba}
B(\a)=\{\b\in \G\,|~0\prec \b \prec \a\}\mbox{ \ \
 for \ }\a\in \G_+.\end{equation}
 The order `` $\succ $" is called {\it dense} if
\begin{equation}\label{dense}
\# B(\a)=\infty \mbox{ \ \ for all \  }\a \in \G_+,\end{equation} it
is  {\it discrete} if
 \begin{equation}\label{discrete}
 B(a)=\emptyset\mbox{ \ \ for some \ }a\in
 \G_+.\end{equation}

Theorem \ref{main-theo2} is an analog of a result in [WS] (where the
Lie algebras $\BB'(\G)$ were considered instead of $\BB(\G)$). Since
the proof is similar to that in [WS], we omit the detail.

Now suppose $\G=\Z$, and from now on, we shall only consider the Lie
algebra
\begin{equation}\label{denote-BB}
\LL=\BB(\Z).\end{equation}

The Lie algebra $\LL$ has a nice realization in the space
$\F[x^{\pm1},t^{\pm1}]\oplus\F C$ by setting
$$L_{i,j}=x^it^{j+1} ~~\mbox{ for any }i, j\in \Z,$$
with the bracket
\begin{equation}\label{realization}
[x^if(t),x^jg(t)]=x^{i+j}(jf'(t)g(t)-if(t)g'(t))+i\d_{i,-j}{\sc\,}\Res{\sc\,}
t^{-1}f(t)g(t),
\end{equation}
for $i,j\in\Z$ and $f(t),g(t)\in\F[t^{\pm1}]$, where the prime
stands for the derivative $\frac{d}{d{\ssc\,}t}$, and
$\Res{\sc\,}f(t)$ stands for the {\it residue} of the Laurent
polynomial $f(t)$, namely the coefficient of $t^{-1}$ in $f(t)$.
Denote by $M(\L)$ the Verma $\BB $-module with highest weight vector
$v_\L$, and by $L(\L)=M(\L)/M'$ the irreducible highest weight
module of weight $\L$, where $M'$ is the maximal proper submodule of
$M(\L)$. Set
\begin{equation}\label{PP}
\AA=\{a\in\BB \,|\,av_\L\in M'\} \mbox{ \ \ and \ \ }\PP=\AA+\BB
_0.\end{equation} Clearly, $\BB _+\subset \AA,$ and $\PP$ is a
subalgebra of $\BB .$

Theorem \ref{main-theo1} will follow from the following proposition.
\begin{prop}\label{prop2} The following conditions are
equivalent:
\begin{itemize}\parskip-4pt

\item[{\rm(1)}] $M(\L)$ is reducible.

\item[{\rm(2)}] $\PP_{-1}\neq \{0\}$.

\item[{\rm(3)}] There exists $0\neq
f(t)=\sum_{j=0}^m a_jt^j\in \F[t]$ for some $m\in\N$ and $a_j\in\F$,
such that $($where $a_{-k}=0$ if $k>0$ or $k<-m)$
\begin{equation}\label{polyss}
\L((t^kf(t))'-a_{-k}C)=0\mbox{ for all }k\in \Z.\end{equation}
\item[{\rm(4)}]
$\Sigma_\L^{(j)}(z):=z\D_\L^{(-j)}(z)-j\D_\L^{(-j-1)}(z)-\frac{c}{j!}z^j$
is a quasipolynomial and satisfies the same nontrivial linear
differential equation with constant coefficients for all $j\in\Z$.
\item[{\rm(5)}]$L(\L)$ is {\it quasifinite}.
\end{itemize}
\end{prop}
\vskip5pt \noindent{\it Proof.~}~We shall follow some arguments in
[S1].

$(1)\Rar(2)$: Assume $M'\ne0$. Similar to the proof in [WS], we can
prove that $\PP$ is a parabolic subalgebra of $\BB $, namely,
$$\PP\supset \BB _0+\BB _+\neq\PP,\mbox{ \ \ and \ }
\PP_{-1}=\PP\cap\BB _{-1}\neq 0.$$

$(2)\Rar(3)$: Let $f(t)$ be the monic ploynomial, called the {\it
characteristic polynomial} of $\PP$, with minimal degree such that
$x^{-1}f(t)\in \PP$ (cf.~[S1, KL]). Set $a=x^{-1}f(t)$. Then from
\begin{equation}\label{t--1}
[t^{-1},x^{-1}tf(t)]=x^{-1}t^{-1}f(t),\ \ \ \
[t^k,x^{-1}f(t)]=-kx^{-1}t^{k-1}f(t)\mbox{ \ \ for \ \ }k\in \Z,
\end{equation} we have $\PP_{-1}=x^{-1}f(t)\F[t^{\pm1}].$ Write
$f(t)=\sum_{j=0}^m a_jt^j$ for some $m\in\N$ and $a_j\in\F$. From
definition (\ref{PP}), we see $M'$ is a proper submodule. We have
$b\cdot av_\L=0$ for all $b\in \BB _+.$ In particular
\begin{equation*}
[xt^k,x^{-1}f(t)]v_\L=\L(-(t^kf(t))'+a_{-k}C)=0~\mbox{ \ for all \
}k\in \Z.\end{equation*}

$(3)\Rar(4)$: Suppose we have (\ref{polyss}). Using
$e^{zt}=\sum_{i=0}^\infty \frac{z^i}{i!}{\sc\,}t^i$ as a generating
series of $\F[t]$, noting that
$f(t)e^{zt}=f(\frac{\ptl}{\ptl\,z})e^{zt}$, by (\ref{polyss}) and
(\ref{realization}), we have (recall that the prime stands for
$\frac{\ptl}{\ptl\,t}$, and the definition of $\D_\L^{(j)}(z)$ in
(\ref{generating-series})),
\begin{eqnarray}
\!\!\!\!\!\!\!\!\!\!0\!\!\!&=\!\!\!&\L((f(t)t^{-j}e^{zt})'-\Res{\sc\,}
t^{-j-1}f(t)e^{zt}C)\nonumber\\\!\!\!\!\!\!\!\!\!\!&=\!\!\!&
\L((f(\frac{\ptl}{\ptl z})t^{-j}e^{zt})')-\Res{\sc\,}
t^{-j-1}f(\frac{\ptl}{\ptl z})e^{zt}c\nonumber
\\\label{generating}\!\!\!\!\!\!\!\!\!\!&=\!\!\!&
\L(-jf(\frac{\ptl}{\ptl z})t^{-j-1}e^{zt}+f(\frac{\ptl}{\ptl
z})zt^{-j}e^{zt})- f(\frac{d}{d
z})\frac{z^j}{j!}c\nonumber\\\!\!\!\!\!\!\!\!\!\!&=\!\!\!&
-jf(\frac{\ptl}{\ptl z})\L(t^{-j-1}e^{zt})+f(\frac{\ptl}{\ptl
z})z\L(t^{-j}e^{zt})- f(\frac{d}{d
z})\frac{z^j}{j!}\nonumber\\\!\!\!\!\!\!\!\!\!\!&=\!\!\!&
-jf(\frac{d}{d z})\D_\L^{(-j-1)}(z)+f(\frac{d}{d
z})z\D_\L^{(-j)}(z)-
f(\frac{d}{d z})\frac{z^j}{j!}c\nonumber\\
\!\!\!\!\!\!\!\!\!\!&=\!\!\!& f(\frac{d}{d z})\Sigma_\L^{(j)}(z)\ \
 \mbox{ \ for }j\in\Z.
\end{eqnarray}
Namely, $\Sigma_\L^{(j)}(z)$ is a quasipolynomial by the statement
after (\ref{generating-series}).

$(4)\Rar(5)$: Assume we have (\ref{generating}).  It suffices to
prove $\PP_{-i}$ has finite codimension in $\BB_{-i}$ for all $i>0$.
For this, we only need to prove $\PP_{-i}\ne0$ for all $i>0$ (then
we can proceed as in (\ref{t--1}) to prove $\PP_{-i}$ has finite
codimension).

 First
from definition (\ref{PP}), $a\in\PP_{-1}$ if and only if $b\cdot
a\cdot v_\L=0$ for all $b\in\BB_1$. Then from (\ref{generating}),
one has \begin{equation}\label{===} x^{-1}t^jf(t)\in\PP_{-1}\mbox{ \
\ for all \ }j\in\Z,\end{equation} this is because $xt^k\cdot
x^{-1}t^jf(t)\cdot v_\L=0$ for all $k\in\Z$, i.e., $\BB_1\cdot
x^{-1}t^jf(t)\cdot v_\L=0$. In particular, from (\ref{===}), one has
$\PP_{-1}\ne0$. Assume that $\PP_{-i}\ne0$ and $0\ne
x^{-i}g(t)\in\PP_{-i}$. Then
$$
x^{-i-1}(-ijt^{j-1}f(t)g(t)-i t^jf'(t)g(t)+t^jf(t)g'(t))=
[x^{-1}t^jf(t),x^{-i}g(t)]\in\PP_{-i-1},
$$
for all $j\in\Z.$ In particular, $\PP_{-i-1}\ne0$.

$(5)\Rar(1)$: It follows from (\ref{infty}).\hfill$\Box$
 \vskip8pt

{\small\ni{\bf Acknowledgement.~}~This work is supported by NSF
grants 10571120, 10471096 of China and Yucai Su's grant ``One
Hundred Talents Program'' from University of Science and Technology
of China. } \vskip8pt

\cl{\bf References}\vs{0pt}

\vskip5pt\small
\parindent=8ex\parskip=2pt\baselineskip=2pt

\re{B} R. Block, On torsion-free abelian groups and Lie algebras,
  {\it Proc. Amer. Math. Soc.} {\bf 9} (1958), 613--620.

\re{DZ} D. Dokovic, K. Zhao, Derivations, isomorphisms and
cohomology of generalized Block algebras, {\it Algebra Colloq.} {\bf
3} (1996), 245--272.


\re{KL} V. Kac, J. Liberati, Unitary quasifinite representations of
$W_\infty$, {\it Lett. Math. Phys.} {\bf 53} (2000), 11--27.

\re{KR} V. Kac, A. Radul, Quasifinite highest weight modules over
the Lie algebra of differential operators on the circle, {\it Comm.
Math. phys.} {\bf 157} (1993), 429--457.

\re{KWY} V.~Kac, W.~Wang, C.H.~Yan, Quasifinite representations of
classical Lie subalgebras of $\WW_{1+\infty}$, {\it Adv.~Math.}
     {\bf139} (1998), 46--140.

\re{LZ} N. Lam, R.B. Zhang, Quasifinite modules for Lie
superalgebras of infinite rank, {\it Tran. Amer. Math. Soc.}
{\bf358} (2005), 403--439.

\re{LW} S. Lin, Y. Wu, Representations of the quantized Weyl algebra
associated to the quantum plane, {\it Algebra Colloquium} {\bf12}
(2005), 715--720.


\re{OZ} J. M. Osborn, K. Zhao, infinite-dimensional Lie algebras of
generalized Block type, {\it Proc. Amer.  Math. Soc.} {\bf 127}
(1999), 1641--1650.

\re{S1} Y. Su, Derivations and structure of the Lie algebras of Xu
type, {\it Manuscripta Math.} {\bf 105} (2001), 483--500.

\re{S2} Y. Su, Structure of Lie superalgebras of Block type related
to locally finite derivations, {\it Comm.  Algebra} {\bf 30} (2002),
3205--3226.

\re{S3} Y. Su, Quasifinite representations of a Lie algebra of Block
type, {\it J. Algebra} {\bf 276} (2004), 117--128.

\re{S4} Y. Su, Quasifinite representations of a family of Lie
algebras of Block type, {\it J. Pure Appl. Algebra} {\bf192} (2004),
293--305.

\re{S5} Y. Su, Classification of quasifinite modules over the Lie
algebras of Weyl type, {\it Adv. Math.} {\bf 174}(2003), 57--68.

\re{S6} Y. Su, Poisson brackets and structure of nongraded
Hamiltonian Lie algebras related to locally-finite derivations, {\it
Canad. J. Math.} {\bf 55} (2003), 856--896.

\re{SX1} Y. Su, X. Xu, Structures of divergence-free Lie algebras,
{\it J. Algebra} {\bf243} (2001), 557--595.

\re{SX2} Y. Su, X. Xu,  Central simple Poisson algebras, {\it
Science in China A} {\bf 47} (2004), 245--263.

 \re{SZ} Y. Su, J.
Zhou, Structure of the Lie algebras with a feature of Block
algebras, {\it Comm. Algebra} {\bf 30} (2002), 3205--3226.


\re{WS} Y. Wu, Y. Su, Highest weight representations of a Lie
algebra of Block type, preprint (arXiv: math.QA/0511733).

\re{X1} X. Xu, Generalizations of Block algebras, {\it Manuscripta
Math}. {\bf 100} (1999), 489--518.

\re{X2} X. Xu, Quadratic conformal superalgebras, {\it J. Algebra}
 {\bf 231} (2000), 1--38.

\re{X3} X. Xu, New generalized simple Lie algebras of Cartan type
over a field with characteristic 0, {\it J. Algebra} {\bf 224}
(2000), 23--58.


\re{ZM} L. Zhu, D. Meng, Structure of degenerate Block algebras,
{\it Algebra Colloq.} {\bf 10} (2003), 53--62.



\end{document}